\documentclass[10pt]{amsart}
\usepackage{amssymb}
\usepackage{epsfig}
\usepackage{amsfonts}
\usepackage{amscd}

\newcommand{\mb}{\overline{\mathcal{M}}}
\newcommand{\cxsm}{(\mathbb{C}^\star)^m}

\newcommand{\tsigma}{\widetilde{\sigma}}
\newcommand{\mbotrei}{\overline{\mathcal{M}}_{0,3}}
\newcommand{\mbotreix}{\overline{\mathcal{M}}_{0,3}(X,d)}

\newcommand{\z}{\mathbb{Z}}

\newcommand{\cx}{\mathbb{C}}

\newtheorem{thm}{Theorem}
\newtheorem{prop}{Proposition}[section]
\newtheorem{cor}[prop]{Corollary}
\newtheorem{lemma}[prop]{Lemma}

\newcommand{\eqqhx}{{QH}^{*}_{T}(X)}

\begin{document}
\title{Positivity in equivariant quantum Schubert Calculus}
\author{Leonardo Constantin Mihalcea }
\begin{abstract} We prove a positivity result in ($T-$)equivariant
quantum cohomology of the homogeneous space $G/P$, generalizing Graham's
positivity in equivariant cohomology.
\end{abstract}
\date{\today} \maketitle
\section{Introduction}
 It is well known that the (integral) cohomology of the
homogeneous space $X=G/P$ (for $G$ connected, semisimple, complex
Lie group and $P$ a parabolic subgroup) satisfies a positivity
property: its structure constants are nonnegative integers equal
to the number of intersection points of three Schubert varieties,
in general position, whose codimensions add up to the dimension of
$X$.

Recently, Graham (\cite{Gr}) has proved a conjecture of Peterson
(\cite{P}), asserting that $H^\star_T(X)$, the $T-$equivariant
cohomology of $X$, where $T \simeq (\cx^\star)^r$ is a maximal
torus in $G$, enjoys a more general positivity property.

Fix $B \supset T$ a Borel subgroup of $G$ and let $W$ and
$W_P$ be the Weyl groups of $G$ and $P$ respectively. Let $\Delta= \{
\alpha_1,...,\alpha_r \}$ be the (positive) simple roots
associated to $W$ and $\Delta_P=\{ \alpha_1,...,\alpha_s \}$ be
those roots in $\Delta$ canonically determined by $P$. The
equivariant cohomology of $X$, denoted $H^\star_T(X)$, is a
$\Lambda=H^\star_T(pt)-$algebra, with a $\Lambda-$basis consisting
of Schubert classes $\sigma(w)^T=[X(w)]_T$. These are determined
by the equivariant Schubert varieties $X(w)_T$ defined with
respect to the Borel group $B$, are indexed by the minimal length
representatives $w$ in $W/W_P$ and have degree $2c(w)$, where
$c(w)$ is the (complex) codimension of $X(w)$ in $X$ (see \S 2
below for more details). The equivariant cohomology of a point is
equal to the polynomial ring $\Lambda=\z[x_1,...,x_r]$ where
$x_i$'s are the negative simple roots in $W$ (see \S 3 below).
Then Graham's Theorem asserts that the
coefficients $c_{u,v}^w$ in $\Lambda$ appearing in the product
\[\sigma(u)^T \cdot \sigma(v)^T = \sum c_{u,v}^w \sigma(w)^T \]
are polynomials in $x_1,...,x_r$ with {\it nonnegative}
coefficients.

The aim of this paper is to show that a similar positivity result holds
in the more general context of $(T-)$equivariant quantum cohomology. This
object was introduced by Givental and Kim (\cite{GK}) for general groups
$G$ and it was initially used to compute presentations of the small
quantum cohomology of (partial) flag manifolds
(\cite{GK,Kim1,Kim2,AS,Kim3}). Related ideas have played a fundamental
role in Mirror Symmetry (\cite{G}).

The equivariant quantum cohomology is a deformation of both the
equivariant and quantum cohomology. It is a graded $\Lambda[q]$-algebra,
where $q=(q_i)$ is indexed by the roots in $\Delta \smallsetminus
\Delta_P$ (or, equivalently, by a basis in $H^2( G/P)$). The multidegree
of $q$ will be given later (\S 4). The basis $\{ \sigma(w)^T \}$ in
equivariant cohomology determines a basis denoted $\{ \sigma(w) \}$ in
its quantum version. If we write $\circ$ for the equivariant quantum
multiplication, the product $\sigma(u) \circ \sigma(v)$ can be expanded
as \[ \sigma(u) \circ \sigma(v) = \sum_d \sum_w q^d c_{u,v}^{w,d}
\sigma(w) \] where $d=(d_i)$ is a multidegree and $q^d$ is the product
$\prod q_i^{d_i}$. The coefficients $c_{u,v}^{w,d}$ are called
equivariant Littlewood-Richardson coefficients (EQLR) and are polynomials
in $\Lambda$ of degree $c(u)+c(v)-c(w)- \sum d_i \cdot \deg(q_i)$. There
is an effective algorithm to compute these coefficients if $X$ is a
Grassmannian, which is a consequence of an equivariant quantum Pieri/Monk
rule (i.e. a multiplication with a divisor class), obtained by the author
in \cite{Mi1}. It is not clear from this algorithm that the EQLR
coefficients enjoy any positivity property. The main result of this paper
is that they do, and, as in the equivariant cohomology, this positivity
holds in the
more general context of homogeneous spaces $G/P$:\\

\noindent { \bf Theorem.} {\it The equivariant quantum
Littlewood-Richardson coefficient $c_{u,v}^{w,d}$ is a polynomial in the
variables $x_1,...,x_r$ with nonnegative coefficients. } \\

\indent The idea of proof is to use a projection formula to reduce
the computation of the EQLR coefficient to a situation where we
can apply Graham's positivity theorem (\cite{Gr}, Theorem 3.2).

{\it Acknowledgments:} I am grateful to my advisor, Prof. W.
Fulton, whose suggesstions greatly influenced the presentation of
this paper. This is part of my thesis.

\section{Cohomology of $G/P$}

In this section we recall the basic facts about the cohomology of
the homogeneous spaces $G/P$. References about this subject can be
found in \cite{Hi} \S 3.3, \cite{FH} especially Ch. 23, \cite{LG},
esp. Chapters 2 and 3 or \cite{Bo}, esp. Ch. 3 and 4. The reader
can also consult the papers \cite{BGG} or \cite{FW}(\S 3).

Unless specified, throughout the paper we use the following
notations (recall from the introduction): $X$ is the homogeneous
space $G/P$, with $G$ complex, connected, semisimple Lie group and
$P$ a parabolic subgroup. $B \subset P$ is a (fixed) Borel
subgroup, $T \simeq (\cx^\star)^r$ the maximal torus in $B$ and
$U$ the unipotent radical of $B$. Recall the Levi decomposition of
$B$ as the semidirect product $T\cdot U$. The set of positive
simple roots $\{ \alpha_1,...,\alpha_r \}$ is denoted by $\Delta$;
$\Delta_P=\{ \alpha_1,...,\alpha_s \}$ denotes the subset of
$\Delta$ canonically determined by the parabolic subgroup $P$. We
write $W=N(T)/T$ for the Weyl group of $G$, where $N(T)$ is the
normalizer of $T$ in $G$. The Weyl group is generated by the
simple reflections $s_{\alpha}$ indexed by elements in $\Delta$.
The parabolic subgroup $P$ determines the subgroup $W_P$ of $W$
generated by all the simple reflections in $\Delta_P$. The length
$l(w)$ of an element $w$ in $W$ is the smallest number of
reflections whose product is $w$; $w_0$ denotes the longest
element in $W$. We write $B^-$ for $w_0 B w_0$, the opposite Borel
subgroup, and $U^-$ for its unipotent radical.

It is a well-known fact that each coset in $W/W_P$ has a unique
representative $w \in W$ of minimal length (see e.g. \cite{LG}, \S
3.5.1); the dual of such $w$, denoted $w^\vee$, is defined to be the
minimal length coset representative of $w_0 w W_P$. Denote by $W^P
\subset W$ the set of minimal length coset representatives for $W/W_P$.
With these notations, $X$ is a smooth variety of dimension equal to the
length of the longest element in $W^P$.

Each $w \in W^P$ determines an element $w P$ in $X$ (in fact, one should
replace $w$ by a representative in $N(T)$, and then consider the
associated coset in $X$. However, since this process is independent of
choices, we denote the result by $w P$). The set $\{ w P \}_{w \in W^P} $
is precisely the set of the $T-$fixed points of $X$ via the
left-multiplication action (\cite{LG}, \S 3.6). For each $w$ in $W^P$ let
$X(w)^o$ be the $U-$orbit $U \cdot wP$ of $wP$ in $X$. This is the affine
Schubert cell, isomorphic (over $\cx$) to the affine space
$\mathbb{A}^{l(w)}$. By the Bruhat decomposition the Schubert cells cover
$X$ with disjoint affines, so their closures $X(w)= cl (X(w)^o)$ (the
Schubert varieties) determine a basis $\{\sigma(w)\}$ for the cohomology
of $G/P$, where $\sigma(w)=[X(w)]$ is of degree $2c(w)=2 (\dim X -
l(w))$. Here $[X(w)]$ denotes the cohomology class in $H^{2c(w)}(X)$
determined by the fundamental class of $X(w)$ via Poincar\'e duality.

Similarly, one can define the {\it opposite} Schubert cell $Y(w)^o \simeq
\mathbb{A}^{l(w)}$ as the $U^--$orbit of $w^\vee P$, determined by the
dual of $w \in W^P$. As before, the opposite Schubert varieties $Y(w)=
cl(Y(w)^o)$ determine a basis $\{\tsigma(w)\}$ for the cohomology of $X$,
where $\tsigma(w)=[Y(w)]$ is in $H^{2c(w)}(X)$. Note that $w_0Y(w)=X(w)$
and since the translations yield the same cohomology class one has that
$\tsigma(w)=\sigma(w)$ in $H^{2c(w)}(X)$.

The bases $\{\sigma(w)\}$ and $\{\tsigma(w)\}$ are dual to each
other in the following sense: if $\pi:X \longrightarrow pt$
denotes the structure morphism of $X$, the cohomology push-forward
$ \pi_\star(\sigma(u) \cdot \tsigma(v))$ is equal to $1 \in
H^0(pt)$ if $u=v^\vee$ and $0$ otherwise\begin{footnote}[1]{ See
Appendix (\S 8 below) for more about this push
forward.}\end{footnote}. This follows from the fact that the
Schubert varieties $X(u)$ and $Y(v)$ intersect properly in $X$,
and, if $c(u)+c(v)= \dim X$, the intersection $X(u) \cap Y(v)$ is
empty unless $u=v^\vee$, when consists of the $T-$fixed point
$uP$, with multiplicity $1$.

Consider the expansion of the product $\sigma(u) \cdot \sigma(v)$ in the
cohomology of $X$:
\[ \sigma(u) \cdot \sigma(v) = \sum_w c_{u,v}^w \sigma(w) .\]
The coefficients $c_{u,v}^w$ are the Littlewood-Richardson coefficients
(LR). They are nonnegative integers, equal to the number of points in the
intersection of $3$ general translates of the Schubert varieties $X(u)$,
$X(v)$ and  $X(w^\vee)$ if $c(u)+c(v)=c(w)$, and $0$ otherwise. The
duality of the bases $\{\sigma(w)\}$ and $\{\tsigma(w)\}$ implies that
the LR coefficients can also be computed as the coefficient of the
fundamental class of a point in $H^0(pt)$:
\[ c_{u,v}^w = \pi_\star(\sigma(u)\cdot \sigma(v) \cdot
\tsigma(w^\vee)).\]

\section{Equivariant cohomology}
\subsection{General facts}
Let $X$ be a complex algebraic variety with a $T\simeq
(\cx^\star)^r-$action. Denote by $p:ET \longrightarrow BT $ the universal
$T-$bundle. $T$ acts on the product $ET \times X$ by $t \cdot
(e,x)=(et^{-1},tx)$. Denote by $X_T$ the quotient $(ET \times X)/T$. The
$T-$equivariant cohomology of $X$, denoted $H^\star_T(X)$, is by
definition the ordinary cohomology $H^\star(X_T)$ of the mixed space
$X_T$. The $X-$bundle $\pi: X_T \longrightarrow BT $ gives the
equivariant cohomology of $X$ a structure of
$H^\star(BT)=H^\star_T(pt)$-algebra.

One can describe the equivariant cohomology of a point as follows (cf.
\cite{Br1}, \S 1): let $\widehat{T}=Hom(T, \cx^\star)$ be the set of
characters of $T$. It is a free abelian group of rank $r$, and the set of
positive simple roots $\Delta=\{ \alpha_1,..., \alpha_r \}$ form a basis.
Each character $\chi \in \widehat{T}$ corresponds canonically to a
$1-$dimensional $T-$module $\cx_\chi$ which determines a line bundle
\[ L(\chi)=ET \times_T \cx_\chi \longrightarrow BT \] over
$BT$, with Chern class $c(\chi)$ ($T$ acts on the right on $ET$, and on
the left on $\cx_{\chi}$). Let $Sym(\widehat{T})$ be the symmetric
algebra of the character group. This is isomorphic to the polynomial ring
$\z[\alpha_1,...,\alpha_r]$ where the complex degree of $\alpha_i$ is
equal to $1$. Then the characteristic homomorphism
\[ c: Sym (\widehat{T}) \longrightarrow H^\star_T(pt) \] sending
$\chi$ to $c(\chi)$ is a ring isomorphism, doubling the degrees.
Recall the notation $\Lambda=H^\star_T(pt)$. In this paper
however, for positivity reasons, we use the {\it negative simple
roots } $ -\alpha_i$, which we denote by $x_i$, as the generators
of $\Lambda$; thus $\Lambda=\z[x_1,...,x_r]$.

There are certain properties of the equivariant cohomology that carry on
from the non-equivariant case. If $X$ is a nonsingular variety (now with
a $T-$action), any $T-$stable subvariety $V$ determines a cohomology
class $[V]_T$ in $H^{2(\dim X - \dim V)}_T(X)$; if $f: X \longrightarrow
Y$ is a $T-$equivariant map of topological spaces, it induces a pull-back
map in cohomology $f^\star:H^i_T(Y) \longrightarrow H^i_T(X)$ for any
integer $i$. In certain situations, for such a $T-$equivariant map, there
is also a Gysin map in cohomology:
\[f_\star^T: H^i_T(X) \longrightarrow H^{i - 2d}_T(Y)\] where $d=
\dim (X) - \dim (Y)$. For the purpose of this paper, we
consider only the situation when $X$ and $Y$ are projective algebraic
varieties and $Y$ is smooth.\begin{footnote}[2]{For more general
situations, such as $X$ or $Y$ noncompact, or being able to find
an ``orientation'' for the map $f$, one can consult e.g.
\cite{FM}, or \cite{F3}, Ch. 19.}
\end{footnote} The definition of this Gysin
map can be found in the Appendix. In what follows we state some of its
properties, in the way they are used in the proof of the positivity
theorem.

{ \it 1. Projection formula:} Let $X_1,X_2$ be projective $T-$varieties,
with $X_2$ smooth, and $f:X_1 \to X_2$ a $T-$equivariant map. Let $a,b$
be equivariant cohomology classes in $H^i_T(X_2)$ and $H^j_T(X_1)$
respectively. Then
\begin{equation}\label{projectionf}f_\star^T(f^\star_T(a)\cdot b)=a \cdot
f_\star^T(b)\end{equation} in $H^{i+j-2d}_T(X_2)$, where $d= \dim X_1 -
\dim X_2$.

{\it 2. Push-forward formula:} Consider the following diagram of
projective $T-$varieties and $T-$equivariant maps:
$$ \begin{CD} X_1 @>{g}>> X_3 \\ @V {f} VV \\  X_2 \end{CD} $$
where $X_2$ and $X_3$ are smooth. Let $V$ be a $T-$invariant subvariety
of $X_3$ of (complex) codimension $c$. Assume that the irreducible
components $V_1,...,V_k$ of $g^{-1}(V)$ have all codimension $c$ in $X_1$
and are $T-$invariant. Let $m_i$ be the (algebraic) multiplicity of $V_i$
in $g^{-1}(V)$. Then
\begin{equation}\label{push} f_\star(g^\star[V]_T) = \sum_i m_i\cdot
a_{i}[f(V_i)]_T \end{equation} in $H^{2c-2d}_T(X_2)$, where $a_{i}$ is a
positive integer equal to the degree of $f_{|V_i}:V_i \longrightarrow
f(V_i)$ or $0$ if $\dim V_i
> \dim f(V_i)$; $d$ denotes again the difference $\dim X_1 - \dim X_2$.\\

Both formulae follow from their non-equivariant counterparts, using the
definition of the equivariant Gysin morphisms, via the finite dimensional
approximations (see \S 8). In the non-equivariant case,
(\ref{projectionf}) is a consequence of the usual projection
formula, while (\ref{push}) follows from the following fact:\\

{\bf Fact:} Let $f: X \to Y$ be a morphism of projective algebraic
varieties, and let $V$ (resp. $W$) be a subvariety of $X$ (resp. $Y$) of
codimension $2c$ (resp. $2k$). Then:

(2') $f_\star[V]=a [f(V)]$ in $H_{2\dim X - 2c}(Y)$ for $[V] \in H_{2\dim
X - 2c}(X)$,
where $a$ is the degree of the map $f_{|V}:V \to f(V)$ or it
is equal to zero if $\dim V$ it is not equal to $\dim f(V)$ (see e.g.
\cite{F1}, Appendix B).
 Here $[V]$ is the homology class
determined by the fundamental class of $V$ in $X$ ($[f(V)]$ is defined similarly).

(2'') Assume $Y$ is smooth and let $[W]$ be the cohomology class in
$H^{2k}(Y)$ determined by $W$. Assume also that $f^{-1}(W)$ has
irreducible components $V_i$ of (the same) codimension $2k$ and algebraic
multiplicity $m_i$. Then
\[ f^\star[W] \cap [X] = \sum_i m_i [V_i] \] in $H_{2\dim X -2k}(X)$. (More
details about this will be given in my thesis \cite{Mi}.) \\

A particular case of an equivariant Gysin map, which will play an
important role in what follows, is the ``integration along the fibres"
$\pi_\star^T:H^i_T(X) \longrightarrow H^{i - 2 \dim (X)}_T(pt)$ induced
by the $T-$equivariant map $\pi: X \longrightarrow pt$. It determines a
$\Lambda-$pairing \[ \langle \cdot, \cdot \rangle_T : H^\star_T(X)
\otimes_\Lambda H^\star_T(X) \longrightarrow \Lambda
\] defined by \[ \langle a,b \rangle_T = \pi_\star^T(a \cup b) \]
More about this pairing will be given in Prop. \ref{duality}
below.

\subsection{Equivariant Schubert calculus on $G/P$}

In this section $X$ denotes the homogeneous space $G/P$. Note that
the Schubert varieties $X(w)$ and $Y(w)$ defined in \S 2 are
$T-$invariant. Since $X$ is smooth, these varieties determine
equivariant cohomology classes $\sigma(w)^T=[X(w)]_T$ and
$\tsigma(w)^T=[Y(w)]_T$ in $H^{2c(w)}_T(X)$. Contrary to the
non-equivariant case, $\sigma(w)^T$ is {\it not} equal to
$\tsigma(w)^T$. \begin{footnote}[3]{In fact, there is an
isomorphism $\overline{\varphi}:H^\star_T(X) \to H^\star_T(X)$
sending $[X(w)]_T$ to $[Y(w)]_T$, induced by the involution
$\varphi:X \to X$ given by $\varphi(x)=w_0\cdot x$. This map is
not $T-$equivariant, but it is equivariant with respect to the map
$T \to T$ defined by $t \to w_0tw_0^{-1}=w_0tw_0$, hence over
$H^\star_T(pt)$ the isomorphism $\overline{\varphi}$ sends
$c(\chi)$ to $c(w_0\chi)$, where $(w_0\chi)(t)=\chi(w_0tw_0)$.}
\end{footnote} Since $\{ \sigma(w) \}_{w \in W^P}$ (resp. $\{ \tsigma(w)
\}_{w \in W^P}$) is a basis for the classical cohomology of $X$, the
Leray-Hirsch Theorem (\cite{Hu} Ch.16), applied to the $X-$bundle $X_T
\to BT$, implies that the set $\{ \sigma(w)^T \}_{w \in W^P}$ (resp. $\{
\tsigma(w)^T \}_{w \in W^P}$) is a basis for the $T-$equivariant
cohomology of $X$. There is an equivariant version of the duality
theorem, described using the $H^\star_T(pt)-$bilinear pairing
\[ \langle \cdot , \cdot \rangle_T : H^\star_T(X)
\otimes_{H^\star_T(pt)} H^\star_T(X) \longrightarrow H^\star_T(pt)\]
given by $x \otimes y \rightarrow \pi_\star^T(x \cdot y)$, where
$\pi_\star^T:H^i_T(X) \longrightarrow H^{i - 2 \dim X}_T(pt)$ is the
integration along the fibers. One has the following result:

\begin{prop}\label{duality}(Equivariant Poincar\'e Duality) The bases
$\{\sigma(w)^T \}$ and $\{ \tsigma(w)^T \}$ are dual to each other i.e.
$\langle \sigma(u)^T, \tsigma(v)^T \rangle_T $ is equal to $1$ if
$u=v^\vee$ and $0$ otherwise.
\end{prop}
\begin{proof} See proof of Lemma 4.2 in \cite{Gr} for $G/B$. The
$G/P$ case is similar. \end{proof}

As in the non-equivariant case, this duality implies a
formula for the {\it equivariant Littlewood-Richardson} (ELR)
coefficients $c_{u,v}^w$ obtained from the expansion
\[ \sigma(u)^T \cdot \sigma(v)^T = \sum_w c_{u,v}^w \sigma(w)^T .
\] They can be computed as follows: \[ c_{u,v}^w =
\pi_\star^T(\sigma(u)^T \cdot \sigma(v)^T \cdot
\tsigma(w^\vee)^T).
\] From either of these descriptions it follows that the ELR
coefficient $c_{u,v}^w$ is a polynomial in
$H^\star_T(pt)=\z[x_1,...,x_r]$ of degree $c(u)+c(v)-c(w)$. Proving a
conjecture of D. Peterson (\cite{P}),  Graham showed that the ELR
coefficients $c_{u,v}^w$ are polynomials in variables $x_1,...,x_r$ with
nonnegative coefficients\begin{footnote}[4]{Graham's result deals with
the case $X=G/B$. The more general situation $X=G/P$ follows from the
fact that the $T-$invariant projection $p:G/B \longrightarrow G/P$
induces an injective map $p^\star_T:H^\star_T(G/P) \longrightarrow
H^\star_T(G/B)$ in equivariant cohomology.} \end{footnote} (cf.
\cite{Gr}). A positive combinatorial formula for these coefficients was
obtained in \cite{KT1} when $X$ is a Grassmannian. The key to that was a
certain recursive formula for the ELR coefficients, which holds in that
case (cf. \cite{MS,O,KT1}). Another recursive formula for any $G/B$ was
obtained in \cite{K}.

\section{Quantum Cohomology of $G/P$}

The (small) quantum cohomology of $X=G/P$ is a graded $\z[q]-$algebra,
having a $\z[q]-$basis consisting of Schubert classes $\sigma(w)$, for
$w$ in $W^P$. Here $q$ stands for the indeterminates sequence $(q_i)$,
indexed by a basis of $H^2(X)$, hence (recall) by the simple roots in
$\Delta \smallsetminus \Delta_P$. The complex degree of $q_i$ is \[ \deg
q_i = \int_{X(s_{\alpha_i})} c_1 (TX) = \pi_\star(\sigma(\alpha_i)\cdot
c_1(TX)) \] where $\pi:X \longrightarrow pt$ is the structure morphism
and $TX$ is the tangent bundle of $X$ (for an explicit computation of
this degree, see e.g. \cite{FW} \S 3). The quantum multiplication,
denoted $\star$, is given by
\[ \sigma(u) \star \sigma(v) = \sum_d \sum_w c_{u,v}^{w,d} q^d
\sigma(w) \] where the first sum is over all sequences of nonnegative
integers $d=(d_i)$ (same number of components as $q$). Recall from the
introduction that $q^d$ stands for the product $\prod q_i^{d_i}$. The
coefficients $c_{u,v}^{w,d}$ are the ($3-$pointed, genus $0$)
Gromov-Witten invariants, equal to the number of rational curves
$f:(\mathbb{P}^1,p_1,p_2,p_3) \longrightarrow X$ of multidegree $d=(d_i)$
from $\mathbb{P}^1$ with three marked points $p_1,p_2,p_3$ to $X$ having
the property that $f(p_1)$ is in $g_1 X(u)$, $f(p_2)$ is in $g_2X(v)$ and
$f(p_3)$ is in $g_3Y(w^\vee)$ for $g_1,g_2,g_3$ general in $G$. This
number is set to be equal to $0$ if $c(u)+c(v)$ is not equal to $c(w)+
\sum \deg q_i \cdot d_i$. We call the coefficient $c_{u,v}^{w,d}$ a
quantum Littlewood-Richardson coefficient.

The fact that such a multiplication gives an associative operation
was proved, using algebro-geometric methods, in \cite{KM} (also
see \cite{FP} and references therein). Computations and properties
of the (small) quantum cohomology algebra were done e.g. in
\cite{W,Be,BCF,Po} for Grassmanians, \cite{C,Ch,FGP,FW,Wo} for
(partial) flag manifolds and in \cite{KTa} for the Lagrangian and
orthogonal Grassmannian. Recently, a simplification of the methods
used to prove results about the quantum cohomology of the
Grassmannian was achieved by Buch (\cite{Bu1}). Similar ideas were
since used in \cite{Bu2,Bu3,BKT1} for more general situations.

We recall an equivalent definition of the coefficient $c_{u,v}^{w,d}$,
which will be generalized in the next section. Let $\mbotreix$ be
Kontsevich' moduli space of stable maps. This is a projective algebraic
variety of dimension $\dim X + \sum d_i \cdot (\deg q_i)$, whose (closed)
points consist of rational maps $f:(C,p_1,p_2,p_3) \longrightarrow X$ of
multidegree $d$, where $C$ is a tree of $\mathbb{P}^1$'s. There are
evaluation maps $ev_i:\mbotreix \longrightarrow X$ which send a point
$(C,p_1,p_2,p_3;f)$ to $f(p_i)$ ($i=1,2,3$), and a contraction map
$\pi:\mbotreix \longrightarrow \mbotrei \simeq pt$ (for details see e.g.
\cite{FP}). Then the quantum LR coefficient is given by
\[ c_{u,v}^{w,d}= \pi_\star\bigl( ev_1^\star(\sigma(u))\cdot
ev_2^\star(\sigma(v)) \cdot ev_3^\star(\widetilde{\sigma}(w^\vee))
\bigr) \] in $H^0(pt)$.

\section{Equivariant quantum cohomology}

The definition of the equivariant quantum cohomology is analogous to the
definition of the quantum cohomology, using now the $T-$equivariant
classes and the equivariant Poincar\'e duality (Prop. \ref{duality}). We
keep the notations from the previous section.

The moduli space of stable maps $\mbotreix$ has a $T-$action given
by:
\[ t \cdot (C,p_1,p_2,p_3;f)
:= (C,p_1,p_2,p_3;\tilde{f}) \] where $\tilde{f}(x):=t \cdot
f(x)$, for $x$ in $C$ and $t$ in $T$. Let $\mbotreix_T$ be its
homotopic quotient, as in \S 3.1. Note that the evaluation maps
$ev_i: \mbotreix \longrightarrow X$ and the contraction map
$\pi:\mbotreix \longrightarrow pt$ are $T-$equivariant. Fix
$u,v,w$ three elements in $W^P$, and $d=(d_i)$ a multidegree.
Following \cite{Kim2} \S 3.1 define the equivariant Gromov-Witten
invariant
\[ c_{u,v}^{w,d} = \pi_\star^T\bigl((ev_1^T)^\star
(\sigma(u)^T) \cdot (ev_2^T)^\star (\sigma(v)^T) \cdot
(ev_3^T)^\star (\tsigma(w^\vee)^T)\bigr) .\] By definition,
$c_{u,v}^{w,d}$ is a homogeneous polynomial of (complex) degree
$c(u)+c(v)-c(w)-\sum d_i \cdot \deg q_i$. Recall that $\Lambda$
denotes $H^\star_T(pt)=\z[x_1,...,x_r]$ (see \S 3.1).

Let $(A, \circ)$ be the graded $\Lambda[q]$-module, where $q=(q_i)$ and
the grading is as before, with a $\Lambda[q]$-basis $\{\sigma(u)\}$
indexed by elements in $W^P$. Define a multiplication, denoted $\circ$,
among the basis elements of $A$ as follows: \[ \sigma(u) \circ \sigma(v)
= \sum_{d}\sum_{w} c_{u,v}^{w,d}q^d \cdot \sigma(w) .
\] The next result was proved by Kim:

\begin{prop}[\cite{Kim2}]\label{definition} $(A,\circ)$ is a
commutative, associative $\Lambda[q]$-algebra with unit. There are
canonical isomorphisms
\begin{enumerate} \item $A/\langle \Lambda^+ \cdot A \rangle
\simeq QH^\star(X)$ as $\z[q]$-algebras. \item $A/\langle q \cdot A
\rangle \simeq H^\star_T(X)$ as $\Lambda$-algebras.
\end{enumerate} sending a basis element $\sigma(u)$ to the corresponding
$\sigma(u)$ in $QH^\star(X)$, respectively to $\sigma(u)^T$ in
$H^\star_T(X)$, where $\Lambda^+$ denotes the ideal of elements in
$\Lambda$ of (strictly) positive degree.
\end{prop}
\begin{proof} The main thing to prove is the associativity. This
was proved in \cite{Kim2} \S 3.3, using a slightly different (but
equivalent) definition of the equivariant Gromov-Witten invariants. A
proof of the equivalence, as well as proofs for the other statements can
be found in \cite{Mi1}, Prop. 3.1, when $X$ is a Grassmannian. The
general case will be given in my thesis.
\end{proof}
\indent It is useful, however, to note that the canonical
isomorphisms in the proposition follow from the the fact that the
equivariant Gromov-Witten invariants specialize to both
equivariant and quantum Littlewood-Richardson coefficients: if
$d=0$, $c_{u,v}^{w,d}$ is the equivariant LR coefficient
$c_{u,v}^w$, while if $c(u)+c(v)=c(w)+\sum d_i \cdot \deg q_i$
(i.e. if $c_{u,v}^{w,d}$ has polynomial degree 0), it is equal to
the quantum LR coefficient $c_{u,v}^{w,d}$.

The algebra $(A,\circ)$ from the previous proposition is the
equivariant quantum cohomology of $X$ and it is denoted by
$\eqqhx$. We call its structure constants, $c_{u,v}^{w,d}$, the
{\it equivariant quantum Littlewood-Richardson} coefficients
(EQLR).

\section{The positivity Theorem}
\subsection{Preliminaries} As usual $X$ denotes the homogeneous
space $G/P$ and $B=T \cdot U$ a  (fixed) Borel subgroup together with its
Levi decomposition (see \S 2). There are two main ingredients in the
proof of the positivity theorem. The first is to consider the space $X
\times X$ with the semidirect product $B'= T \cdot (U^- \times U^-)$
acting by
\[t\cdot (u_1,u_2)(\bar{g}_1,\bar{g}_2)=(tu_1\bar{g}_1,tu_2\bar{g}_2)\]
($(\bar{g}_1,\bar{g}_2) \in X \times X$), then to prove that any
subvariety $V$ of $X \times X$ satisfies a positivity property.
This is where we use Graham's results. The second ingredient is to
reduce the computation of an EQLR coefficient to a computation
taking place on $X \times X$, and then use the positivity property
there. This is done in \S 6.2 , using a projection formula. We
begin by stating a special case of the key result in \cite{Gr}:

\begin{prop}[see \cite{Gr} Thm. 3.2]\label{positivity} Let $X \times X$ with
the $B'$ action defined above. Let $\beta_1,...,\beta_d$ in the character
group of $T$ denote the weights of the adjoint action of $T$ on the Lie
algebra $Lie(U^- \times U^-)$ of the unipotent radical of $B'$. Let $V$
be a $T-$stable subvariety of $X \times X$. Then there exist $B'-$stable
subvarieties $D_1,...,D_t$ of $X \times X$ such that in $H^\star_{T}(X
\times X)$,
\[ [V]_{T}= \sum f_i [D_i]_{T} \] where each $f_i \in H^\star_{T}(pt)$
can be written as a linear combination of monomials in
$\beta_1,...,\beta_d$ with nonnegative integer coefficients.
\end{prop}

\noindent {\it Remark:} Graham's result deals with the more general
situation of a variety $Y$, possibly singular, a connected, solvable
group $B'=T'U'$ acting on it and a $T'-$stable subvariety $V$ of $Y$. $V$
determines only an equivariant {\it homology} class $[V]_{T'} \in
H^{T'}_{2 \dim V}(Y)$, therefore his positivity result takes place in
$H^{T'}_\star(Y)$. If $Y$ is smooth the equivariant homology and
cohomology are identified via the equivariant Poincar\'e duality
(\cite{Br2} \S 1), and one recovers Prop. \ref{positivity} (for $Y= X
\times X$).

\begin{cor}\label{mypositivity} Let $X \times X$ endowed with the previous
$B'-$action, and $V$ a $T-$stable subvariety. Then $[V]_T$ can be
written uniquely as
\[ [V]_T = \sum f_i [Y(u) \times Y(v)]_T \] in $H^\star_T(X \times X)$
where each $f_i$ is a polynomial in $x_1,...,x_r$ with nonnegative
coefficients.
\end{cor}

\begin{proof} Note first that the weights of $T$ acting
on $Lie(U')$ ($=Lie(U^-) \times Lie(U^-)$) are the same as the weights of
$T$ acting on $Lie(U^-)$, which are the negative roots of $G$. To finish
the proof, it is enough to show that the $B'-$stable subvarieties of $X
\times X$ are precisely the products of Schubert varieties $X(u) \times
X(v)$, and that they determine a basis for the equivariant cohomology of
$X \times X$.

To do that, recall that $X$ has finitely many $U^--$orbits, the
affine Schubert cells $Y(w)^o$. These orbits are $B^--$stable
(they are $B^--$orbits of the $T-$fixed
points)\begin{footnote}[5]{The fact that $X$ has finitely many
orbits implies that any $U^--$orbit is $B^--$stable holds in a
more general context, see \cite{Gr}, Lemma 3.3.}\end{footnote},
and they cover $X$ with disjoint affines. Then the unipotent
radical $U'=U^- \times U^-$ of $B'$ acts on $X \times X$ with
finitely many orbits $Y(u)^o \times Y(v)^o$ ($u,v \in W^P$), and
these orbits cover $X \times X$ by disjoint affines. Their
closures are products of Schubert varieties $Y(u) \times Y(v)$,
and determine a basis $[Y(u) \times Y(v)]_T$ for the equivariant
cohomology $H^\star_T(X \times X)$. Moreover, since any $U'-$orbit
is $B'-$stable, it follows that any irreducible $B'-$stable
variety must be one of $Y(u) \times Y(v)$. \end{proof}

We state next a result which is based on Kleiman's transversality theorem
and on a lemma about translates of Schubert varieties (cf. \cite{FW} \S
7):

\begin{lemma}[\cite{FW}, Lemma 7.2]\label{Kleiman} Let $Z$ be an irreducible
$G-$variety and let $F:Z \longrightarrow X \times X$ be a $G-$equivariant
morphism, where $G$ acts diagonally on $X \times X$. Then, for any $u$
and $v$ in $W^P$, the subscheme $F^{-1}(X(u) \times Y(v))$ is reduced,
locally irreducible, of codimension $c(u)+c(v)$.
\end{lemma}

Given $v \in W^P$, we apply Lemma \ref{Kleiman} to $Z=\mbotreix$ (with
the $G-$action induced from $X$, as in \S 5), $F=(ev_3,ev_3)$ (recall
that $ev_3$ is the evaluation map, cf. \S 4) and $u=w_0$ (the longest
element in $W^P$). $\mbotreix$ is irreducible by \cite{KP,T} and $F$ is
clearly $G-$equivariant. Then $ev_3^{-1}(Y(v))=F^{-1}(X \times Y(v))$ is
reduced, locally irreducible, of codimension $c(v)$. Next result shows
slightly more:

\begin{lemma}\label{inverse} The inverse image $ev_3^{-1}(Y(v))$ of $Y(v)$ is a
disjoint union of $T-$stable, reduced, irreducible components,
each of them with codimension equal to $c(v)$, the codimension of
$Y(v)$ in $X$.
\end{lemma}

\begin{proof} It remains to show that each
component $V_i$ of $ev_3^{-1}(Y(v))$ (which is also a connected
component) is $T-$stable. Since the whole preimage is $T-$stable,
it follows that the $T-$action sends one connected component to
another. But the identity $id$ in $T$ fixes all the $V_i$'s so
each component must be $T-$stable.\end{proof}

\subsection{Proof of the Theorem} The idea of proof is to use the projection
formula (\ref {projectionf}) and the push-forward formula (\ref
{push}) from \S 3.1 to reduce the computation of the EQLR
coefficients from an (equivariant) intersection problem on a
moduli space to an (equivariant) intersection problem on the
product $X \times X$, where we apply Corollary \ref{mypositivity}.
If $Y$ is a variety, denote by $\pi_Y$ the structure morphism
$\pi: Y \longrightarrow pt$. Denote by $\mb$ the moduli space
$\mbotreix$ and by $\pi_{\mb}$ the morphism $\pi_{\mbotreix}$.

\begin{thm}\label{positivity Thm} Let $u,v,w$ in $W^P$ and $d$ a multidegree.
Then the equivariant quantum Littlewood-Richardson coefficient
$c_{u,v}^{w,d}$ is a polynomial in the (negative) simple roots
$x_1,...,x_r$ with nonnegative coefficients. \end{thm}

\begin{proof} Let $F:\mb \longrightarrow X \times X$ be
$(ev_1,ev_2)$. Clearly, $F$ is $T-$equivariant and proper. The main
point of the proof is the following Claim: \\

\noindent {\it Claim: The EQLR coefficient $c_{u,v}^{w,d}$ is equal to
\[\sum a_i (\pi_{X \times X})_\star^T\bigl([X(u)\times X(v)]_T \cdot
[F(V_i)]_T\bigr) \]in $H^\star_T(pt)$, where the sum is over the
components $V_i$ of $ev_3^{-1}(Y(w^\vee))$ and $a_i$ is the degree
of the map $F_{|V_i}:V_i \to F(V_i)$ or $0$ if $\dim V_i > \dim F(V_i)$.}

We postpone the proof of the Claim until later, and we show first how
does the Claim imply the theorem. Clearly, it is enough to show that each

\begin{eqnarray}\label{*}(\pi_{X \times X})_\star^T\bigl([X(u)\times X(v)]_T \cdot
[F(V_i)]_T\bigr)\end{eqnarray}

\noindent is nonnegative in the sense of the theorem. By Corollary
\ref{mypositivity}, each $[F(V_i)]_T$ can be written as a combination
$\sum f^i_{j} [Y(u_j) \times Y(v_j)]_T$, with $u_j,v_j$ in $W^P$, where
each $f^i_j$ is a nonnegative polynomial in variables $x_1,...,x_r$.
Projection formula (\ref{projectionf}) implies that the expression
(\ref{*}) can be written as a combination
\[\sum P^i_j (\pi_{X \times X})_\star^T\bigl([X(u) \times X(v)]_T \cdot [Y(u_j) \times
Y(v_j)]_T\bigr)\] with each of $P^i_j$ a nonnegative sum of $f^i_j$'s,
hence also nonnegative. Write $pr_i:X \times X \longrightarrow X$ for the
projection to the $i-$th component of $X$ ($i=1,2$). Then
\[(\pi_{X \times X})_\star^T\bigl([X(u) \times X(v)]_T \cdot [Y(u_j)
\times Y(v_j)]_T\bigr)\] can be written as
\begin{eqnarray*} (\pi_{X \times X})_\star^T\bigl([X(u)
\times X(v)]_T \cdot [Y(u_j) \times Y(v_j)]_T\bigr)& =& \\ (\pi_{X \times
X})_\star^T\bigl((pr_1^T)^\star(\sigma(u)^T) \cdot
(pr_2^T)^\star(\sigma(v)^T) \cdot [Y(u_j) \times Y(v_j)]_T\bigr) & = &
\\ (\pi_X)_\star^T\bigl(\sigma(u)^T \cdot [Y(u_j)]_T\bigr) \cdot
(\pi_X)_\star^T\bigl(\sigma(v)^T \cdot [Y(v_j)]_T\bigr)  & = &
\delta_{u^\vee,u_j} \cdot \delta_{v^\vee,v_j}
\end{eqnarray*}

The first equality follows from the fact that $[X(u) \times X(v)]_T =
(pr_1^T)^\star(\sigma(u)^T) \cdot (pr_2^T)^\star(\sigma(v)^T)$, the
second equality is Lemma 2.3 from \cite{Gr} and the third equality
follows from the equivariant Poincar\'e duality (Prop. \ref{duality}). This
proves the positivity theorem.
\end{proof}

\begin{proof}[Proof of the Claim] Recall (\S 5) that $c_{u,v}^{w,d}$
is defined by
\[ c_{u,v}^{w,d} = (\pi_{\mb})_\star^T\bigl((ev_1^T)^\star
(\sigma(u)^T) \cdot (ev_2^T)^\star (\sigma(v)^T) \cdot (ev_3^T)^\star
(\sigma(w^\vee)^T)\bigr) \] Consider the composite
$$ \begin{CD} \mb_T @>{(F^T,ev_3^T)}>>(X \times X)_T \times
 X_T \simeq (X_T \times_{BT} X_T) \times X_T \end{CD} $$
Then
\begin{eqnarray*} (ev_1^T)^\star (\sigma(u)^T) \cdot
(ev_2^T)^\star (\sigma(v)^T) \cdot (ev_3^T)^\star (\tsigma(w^\vee)^T) &=&
(F^T)^\star\bigl([X(u) \times X(v)]_T\bigr) \cdot (ev_3^T)^\star
(\tsigma(w^\vee)^T) \end{eqnarray*} But the composite

$$ \begin{CD} \mb@>{F}>>X\times X @>{\pi_{X \times
X}}>>pt \end{CD} $$ is equal to $\pi_{\mb}$ , therefore
$c_{u,v}^{w,d}$ is equal to

\[c_{u,v}^{w,d} = (\pi_{ X \times X})_\star^T F^T_\star\bigl(
(F^{T})^\star([X(u) \times X(v)]_T) \cdot (ev_3^T)^\star
(\tsigma(w^\vee)^T)\bigr)\] By projection formula (\ref{projectionf})
\[ F^T_\star\bigl((F^T)^\star([X(u) \times X(v)]_T) \cdot (ev_3^T)^\star
(\tsigma(w^\vee)^T)\bigr) = [X(u) \times X(v)]_T \cdot F^T_\star\bigl(
(ev_3^T)^\star (\tsigma(w^\vee)^T)\bigr) \]

By Lemma 6.3, $ev_3^{-1}(X(w^\vee))$ is a disjoint union of
$T-$stable, reduced, irreducible components $V_i$. Then, by the
push-forward formula (\ref{push}),
\[ F^T_\star( (ev_3^T)^\star [Y(w^\vee)]_T) = \sum a_{i} [F(V_i)]_T \]
where $a_{i}$ are as in the Claim. To conclude, we have obtained that
\[ c_{u,v}^{w,d} = \sum a_{i} (\pi_{X \times X})_\star^T\bigl([X(u)\times X(v)]_T \cdot
[F(V_i)]_T\bigr) \] which proves the Claim. \end{proof}
\section{Positivity for a different torus action}

Let $X$ be the variety of partial flags $F_{m_1} \subset ... \subset
F_{m_k}$ in $\cx^m$, seen as the quotient $PGL(m)/P$ ($P$ parabolic).
There is a $T'=\cxsm$-action on $X$ induced by the one of the maximal
torus $T= \cxsm/\cx^\star$ in $PGL(m)$, via the quotient morphism
$\varphi:\cxsm \to \cxsm/\cx^\star$. Hence for $t' \in T'$ and $x \in X$
\[ t' \cdot x = \varphi(t') x .\]

The goal of this section is to interpret the positivity result
using this action, which is the one used in \cite{Mi1} when $X$
was a Grassmannian. Besides using a different action in {\it loc.
cit.}, we have also used different (but standard) generators for
$H^\star_{T'}(pt)$. We recall their definition. Let $p: ET'_n
\longrightarrow BT'_n$ be the family of finite-dimensional
$T'-$bundles
\[ p : (\cx^{n+1}\smallsetminus 0)^m \longrightarrow \prod_{i=1}^m
\mathbb{P}^n \] approximating the universal $T'-$bundle $ET' \to BT'$
(\cite{Hu} \S 4.11 or Ch. 7, \cite{Br1,Br2},\cite{EG} \S 3.1). If $Y$ is
a $T'-$space, the equivariant cohomology $H^\star_{T'}(Y)$ can be
computed as a limit of the ordinary cohomology of the finite dimensional
approximations: $H^i_{T'}(Y) = \lim_{ n \rightarrow \infty}
H^i(Y_{T',n})$, where $Y_{T',n} = ET'_n \times_{T'} Y$. Then the
equivariant cohomology of a point is
$H^\star_{T'}(pt)=\lim_{n \to \infty}H^\star(BT'_n)=\z[T_1,...,T_m]$, where $T_j$ is the
first Chern class of the Serre bundle $\mathcal{O}(1)$ on the $j-$th
factor of $\prod_{i=1}^m \mathbb{P}^n$.

The EQLR coefficients in $T'-$equivariant quantum cohomology are obtained
as the images of those in $T-$equivariant quantum cohomology via the
morphism \[\overline{\varphi}: H^i_{T}(pt) \to H^i_{T'}(pt)\] induced by
$\varphi$. Using the identification of the $H^\star_T(pt)$ (resp.
$H^\star_{T'}(pt)$) with the character group $\widehat{T}=Hom(T,
\cx^\star)$ (resp. with $\widehat{T'}=Hom(T', \cx^\star)$) (cf. \S 3.1),
one sees that $\overline{\varphi}$ sends the first Chern class of the
line bundle $L(x_j)$ associated to the negative simple root $x_j$ in
$H^2_T(pt)$ (so $x_j$ is the character
$x_j\bigl((z_1,...,z_m)\cx^\star\bigr)=\frac{z_{j+1}}{z_{j}}$) to the
first Chern class of the line bundle corresponding to the canonical lift
of this character to $\widehat{T'}$. It turns out that the latter Chern
class is equal to $T_{j} - T_{j+1}$, defined in the previous paragraph.
Then the positivity theorem implies:
\begin{cor} The EQLR coefficient $c_{u,v}^{w,d}$ in the
$T'-$equivariant quantum cohomology of $X$ is a polynomial in the
variables $T_1-T_2,...,T_{m-1}-T_m$ with nonnegative coefficients.
\end{cor}

\section{Appendix - Equivariant Gysin morphisms} The aim of this Appendix
is to define the equivariant Gysin maps used in \S 3. Let $f:X
\longrightarrow Y$ be a morphism of projective varieties, with $Y$
smooth. Let $d= \dim (X) - \dim (Y) $ (complex dimensions). Define a
Gysin map $f_\star:H^i(X) \longrightarrow H^{i - 2d} (Y)$ by the
composite
$$ \begin{CD} H^i(X) @> { \cap [X] } >> H_{2 \dim (X) - i}
(X) @> {f_\star} >> H_{ 2 \dim (X) - i} (Y) \simeq H^{i - 2d}(Y)
\end{CD} $$ where $[X]$ is the fundamental class of $X$ in the singular homology
group $H_{2 \dim X}(X)$, and the middle $f_\star$ is the singular
homology push-forward (if $X$ or $Y$ were not compact, one should use
Borel-Moore homology). The last isomorphism is given by Poincar\'e
duality. We need the following property of the Gysin map:

\begin{lemma}\label{fiber square} Consider the following fiber square of projective
varieties:
$$ \begin{CD} X' @>{i}>> X \\ @V {f'}VV @V {f} VV \\ Y' @>{j}>>Y
\end{CD} $$ where $Y,Y'$ are smooth and $i,j$ are regular
embeddings of the same (complex) codimension $c$. Then $f'_\star i^\star
= j^\star f_\star$ as maps $H^i(X) \to H^{i-2d}(Y')$.
\end{lemma}

\begin{proof} The proof is given in my thesis \cite{Mi} (one could also see \cite{FM}).
\end{proof}

Assume the map $f:X \to Y$ (with $Y$ smooth) is $T-$equivariant.
Then it determines a Gysin map of the cohomology of the
finite-dimensional approximations (see \S 7)
$f_{\star,n}:H^i(X_{T,n}) \longrightarrow H^{i-2d} (Y_{T,n})$.
Define the equivariant Gysin map $f_\star^T:H^i_T(X)
\longrightarrow H^{i-2d}_T(Y)$ as the unique map that makes the
following diagram commute:
$$ \begin{CD} H^i(X_{T,n}) @<{res}<<H^i_T(X) \\
@V{f_{\star,n}}VV @V{f_\star^T}VV\\ H^{i- 2d}(Y_{T,n})
@<{res}<<H^{i-2d}_T(Y) \end{CD} $$for every integer $n$. The
horizontal maps $res$ are the cohomology pull-backs induced by the
inclusions $X_{T,n} \to X_T$ (resp. $Y_{T,n} \to Y_T$). The
existence and uniqueness of $f_\star^T$ follow from the fact that
the equivariant cohomology can be computed by passing to the limit
on the ordinary cohomology of the finite dimensional
approximations. The only thing one has to check is that the the
ordinary Gysin maps are compatible to each other, in the sense
that the diagram
\begin{equation*}
\begin{CD}
H^i(X_{T, n_1}) @< {res_{}}<<H^i(X_{T, n_2})\\@V{f_{\star,n_1}}VV @V{f_{\star,n_2}}VV \\
H^{i-2d}(Y_{T, n_1}) @<{res_{}}<< H^{i-2d}(Y_{T, n_2})
\end{CD}
\end{equation*} is commutative for any integers $n_1 < n_2$. This follows by applying Lemma
\ref{fiber square} to the fiber square
\begin{equation*} \begin{CD}
X_{T, n_1} @> {i}>>X_{T, n_2}\\@V{f_{n_1}}VV @V{f_{n_2}}VV \\
Y_{T, n_1} @>{j}>> Y_{T, n_2} \end{CD} \end{equation*}

\end{document}